\documentclass[11pt]{amsart}
\usepackage{amscd,amssymb,amsmath}
\usepackage[all]{xy}
\theoremstyle{definition}

\theoremstyle{remark}

\numberwithin{equation}{section}

\newcommand{\ga}{\Gamma}

\newcommand{\bea} {\begin{eqnarray*}}
\newcommand{\beq} {\begin{equation}}
\newcommand{\bey} {\begin{eqnarray}}
\newcommand{\eea} {\end{eqnarray*}}
\newcommand{\eeq} {\end{equation}}
\newcommand{\eey} {\end{eqnarray}}
\newcommand{\til}[1]{\widetilde{#1}}
\newcommand{\ovl}{\overline}

\def\Z{{\mathbb Z}}

\bibliography{ref}
\begin{document}
\title[Twisted conjugacy classes]
{Twisted conjugacy classes in Symplectic groups, Mapping class
 groups and Braid groups\\
(including an Appendix written  with Francois Dahmani)}

\author{Alexander Fel'shtyn}
\address{ Instytut Matematyki, Uniwersytet Szczecinski,
ul. Wielkopolska 15, 70-451 Szczecin, Poland
and Boise State University, 1910
University Drive, Boise, Idaho, 83725-155, USA }
\email{felshtyn@diamond.boisestate.edu, felshtyn@mpim-bonn.mpg.de}
\author {Daciberg L. Gon\c{c}alves}
\address{Dept. de Matem\'atica - IME - USP, Caixa Postal 66.281 - CEP
05311-970,
S\~ao Paulo - SP, Brasil}
\email{dlgoncal@ime.usp.br}
\address{ Francois Dahmani, Laboratoire Emile Picard, Umr CNRS 5580, Universit\'e Paul
 Sabatier,
31062 Toulouse Cedex 4, France}
\email{francois.dahmani@math.ups-tlse.fr}
\thanks{This work was initiated during our  visit  to the   University of  British Columbia, Vancouver in April 2007 and the visit
 of the first  author to the Universit\'e Paul Sabatier,  Toulouse in May-June 2007  and was completed during conferences in Bendlewo and Warsaw in July-August 2007. }.

\begin{abstract}
 We prove  that the symplectic group  $Sp(2n,\mathbb Z)$
and the mapping class group  $Mod_{S}$ of a compact  surface $S$ satisfy  the $R_{\infty}$ property. We also show that $B_n(S)$, the full braid group on $n$-strings of a  surface  $S$, satisfies the $R_{\infty}$ property in the  cases where $S$ is either the compact  disk $D$, or the sphere $S^2$. This means that  for any automorphism $\phi$ of $G$, where $G$ is one of the  above groups,   the number of twisted $\phi$-conjugacy classes  is infinite.
\end{abstract}

\date{\today}
\keywords{Reidemeister number, twisted conjugacy classes, Braids group,
 Mapping class group, Symplectic group}
\subjclass[2000]{20E45;37C25; 55M20}
\maketitle

\tableofcontents

\section{Introduction }
Let $\phi:G\to G$ be an automorphism of a group $G$. A class of
 equivalence relation defined by
$x\sim gx\phi(g^{-1})$ is called the \emph{Reidemeister class },
 $\phi$-\emph{conjugacy class}, or \emph{twisted conjugacy class of}
 $\phi$. The
number $R(\phi)$ of Reidemeister classes is called the
 \emph{Reidemeister number}
of $\phi$.
The interest in twisted conjugacy relations has its origins, in
 particular,
in the Nielsen-Reidemeister fixed point theory (see, e.g.
 \cite{reid:re,j, FelshB}),
in Selberg theory (see, eg. \cite{Shokra,Arthur}), and in  Algebraic
 Geometry
(see, e.g. \cite{Groth}).

A current significant  problem in this area is to obtain
a twisted analogue of the  Burnside-Frobenius theorem
\cite{FelHill,FelshB,FelTro,FelTroVer,ncrmkwb,polyc,FelTroObzo},
that is, to show the
coincidence of
the Reidemeister number of $\phi$ and the number of fixed points of the
induced homeomorphism of an appropriate dual object.
One step in this process is to describe the class of groups $G$
 for which  $R(\phi)=\infty$ for
any automorphism $\phi:G\to G$.

The work of discovering which groups belong to this  class of
 groups was begun by Fel'shtyn and Hill in \cite{FelHill}.
 It was later  shown by various authors that the following groups belong
 to this class:
 (1) non-elementary Gromov hyperbolic groups \cite{FelPOMI,ll}; (2)
Baumslag-Solitar groups $BS(m,n) = \langle a,b | ba^mb^{-1} = a^n
 \rangle$
except for $BS(1,1)$ \cite{FelGon};
 (3) generalized Baumslag-Solitar groups, that is, finitely generated
 groups
which act on a tree with all edge and vertex stabilizers infinite
 cyclic
\cite{LevittBaums}; (4)
lamplighter groups $\mathbb Z_n \wr \mathbb Z$ if and only if  $2|n$ or $3|n$
 \cite{gowon1}; (5)
the solvable generalization $\ga$ of $BS(1,n)$ given by the short exact
 sequence
$1 \rightarrow \mathbb Z[\frac{1}{n}] \rightarrow \ga \rightarrow
 \mathbb Z^k \rightarrow 1$
 as well as any group quasi-isometric to $\ga$ \cite{TabWong},
such groups  are quasi-isometric to $BS(1,n)$ \cite{TabWong2} ( note however that the class of groups
for which $R(\phi)=\infty$ for
any automorphism $\phi$ is not closed under quasi-isometry);  (6) saturated weakly branch groups, including the Grigorchuk group and the Gupta-Sidki group
 \cite{FelLeoTro};
(7)The  R. Thompson group F \cite{bfg}.

The paper \cite{TabWong} suggests a terminology for this property,
 which we would like to follow.  Namely, a group $G$ has \emph{property }
 $R_\infty$ if all of its automorphisms $\phi$
satisfy  $R(\phi)=\infty$.

For the immediate consequences of the $R_\infty$ property in
 topological fixed point
theory see, e.g., \cite{TabWong2}.

In the present   paper we prove  that the symplectic group  $Sp(2n,\mathbb Z)$
and the mapping class group  $Mod_{S}$ of a compact  surface $S$ have  the $R_{\infty}$ property. We also show that $B_n(S)$, the full braid group on $n$-strings of a compact  surface $S$, satisfy the $R_{\infty}$ property in the  cases where
 $S$ is either the compact  disk $D$, or the sphere $S^2$.
The results of the present paper indicate  that the further study
of Reidemeister theory for these groups  should  go
along the lines similar to those  of  the infinite case. On the other
hand, these  results  reduces  the class of groups for which the
twisted Burnside-Frobenius conjecture
 \cite{FelHill,FelTro,FelTroVer,ncrmkwb,polyc,FelTroObzo}
has yet to be verified.

The paper is organized into 6  sections. In section 2 we describe a very
 naive procedure to decide  whether  a group has the $R_{\infty}$ property.
 Also we recall  some   known relations  between mapping class groups  and braid
 groups  which will be used later.
In section 3  we show that the symplectic group
 has the $R_{\infty}$ property. In section 4  we give   algebraic
 proof of the $R_{\infty}$ property for the Mapping class group of  the  closed surfaces. There we  use  the short exact sequence
$$1\to I_n \to Mod_S \to Sp(2n,\mathbb Z)\to 1,$$
 where $I_n$ is the Torelli group.
 In section 5  we show that,  with a  few exceptions, the braids groups of
 the disk  $D$ and the sphere $S^2$  have the $R_{\infty}$ property.

In Appendix, written with Francois Dahmani, we use   geometric  methods to  show that  the mapping class group  of the compact surfaces  has  the $R_{\infty}$ property  with few obvious exceptions.

{\bf Acknowledgments} The first  author would like to thank
V. Guirardel,  J. Guaschi,  M. Kapovich and J.D. McCarthy    for stimulating
 discussions and comments. He also thanks   the Universit\'e Paul Sabatier, Toulouse  for its kind  hospitality and support while the  part of this work has been completed.

\vskip 1.0cm

\section{ Preliminaries}

Consider a group extension respecting the  homomorphism $\phi$:

\begin{equation}\label{extens}
\xymatrix{
0\ar[r]&
H \ar[r]^i \ar[d]_{\phi'}&  G\ar[r]^p \ar[d]^{\phi} & G/H \ar[d]^{\ovl{\phi}}
\ar[r]&0\\
0\ar[r]&H\ar[r]^i & G\ar[r]^p &G/H\ar[r]& 0,}
\end{equation}
where $H$ is a normal subgroup of $G$.
First,  let us remark  that the Reidemeister classes of $\phi$ in $G$
are mapped epimorphically on classes of $\ovl\phi$ in $G/H$. Indeed,
\begin{equation}
p(\til g) p(g) \ovl\phi(p(\til g^{-1}))= p (\til g g \phi(\til g^{-1}).
\end{equation}

Moreover, if  $R(\phi)<\infty$, then the previous remark implies
$R(\ovl\phi)<\infty$.

 We next give    a criterion for a group  to have
 the $R_{\infty}$ property.

{\bf Lemma 2.1} Suppose that a group $G$ has an infinite number of
 conjugacy classes. Then  any inner automorphism also has infinite
  Reidemeister number.
Moreover, if $\{\phi_j\}_{j\in Out(G)}$ is a subset of $Aut(G)$ which
 contains single   representatives
 for each  coset $Aut(G)/Inn(G)$ and $R(\phi_j)=\infty$ for all $j$,
 then $G$ has the
 $R_{\infty}$ property.
\begin{proof} The  first part of the statement was proved  in
 \cite{bfg}. The second part follows easily using a similar argument which we now recall.
 Let $\phi$
 be an automorphism and  $\theta$ an element of the group. Then we have
 two automorphisms, namely $\phi$ and the composite  of $\phi$  with the
 inner automorphism which is conjugation by $\theta$, which we denote
 by $\theta\circ\phi$.  We  claim that multiplication by
$\theta^{-1}$ on the right provides a bijection between the set of
 Reidemeister classes of $\phi$ and  those  of
 $\theta\circ\phi$. Consider two elements where the first is denoted by
  $\alpha$ and the second is of the form $\beta\alpha\phi(\beta)^{-1}$
 for
 some $\beta\in G$. We  claim that the two elements   $\alpha\theta^{-1}$
 and the
 $\beta\alpha\phi(\beta)^{-1}\theta^{-1}$ are in the same
  $\theta\circ\phi$ Reidemeister classes.
To show this,  we write 
 class of  $\alpha\theta^{-1}$  as
 $\beta\alpha\theta^{-1}
 ((\theta\circ\phi(\beta))^{-1})=\beta\alpha\theta^{-1}
\theta\phi(\beta)^{-1}\theta^{-1}=\beta\alpha\phi(\beta)^{-1}\theta^{-1}$. Thus  $\alpha\theta^{-1}$ and
 $\beta\alpha\phi(\beta)^{-1}\theta^{-1}$ are in the same Reidemeister
 class of $\theta\circ\phi$.
 Similarly,  multiplication by $\theta$ on the right provides a bijection
 between the set of Reidemeister classes of  $\theta\circ\phi$ and  those of
 $\phi$.
 One correspondence is the inverse of the other and the result follows.
\end{proof}

{\bf Lemma 2.2} Let $\phi\in Aut(G)$ and $Q_i$ be  an infinite family of
  quotients of $G$,  $$1\to K_i \to G\to Q_i \to 1,$$
\noindent such that $\phi(K_i)\subset K_i$. If the  sequence of numbers
 $R(\phi_i)$ is unbounded, where
$\phi_i$ is the induced homomorphism on the quotient $Q_i$, then
 $R(\phi)=\infty$.
\begin{proof}
 Since  $p:G \to Q_i$ is surjective, we have
 $R(\phi)\geq R(\phi_i)$ . If $R(\phi)$ is finite then  the
 sequence
$R(\phi_i)$ is bounded. Hence  the result follows.
\end{proof}

Despite the fact that the Lemma 2.2  is quite obvious, it is
 usefull for its  ability to find such quotients
 and   to estimate Reidemeister number on the quotient. We will apply
 Lemma 2.2  in   a situation  where the   quotients are  finite groups.

 We will now  state a result   which relates
 braid groups and mapping class groups. We will use this relation  to
 study the braid group of the  sphere. Following section 2
 of \cite{Scott}, let $S$ be a 2-manifold
with $n$ distinguished points in its interior. Let  $H(S,n)$ denotes
 the space of homeomorphisms of $S$, which fix pointwise the  $n$  distinguished points. If $n=0$,
 we set  $H(S)=H(S,0)$. We define
$G(S,n)=\pi_0(H(S,n))$ and,  if $n=0$, we set  $G(S)=G(S,0)$.

If $S$ is the sphere $S^2$ we have (see the bottom of the page 615 in
 \cite{Scott}):

{\bf Theorem 2.3} If $S$ is  $S^2$ then we have the  short exact
 sequence

$$1\to \mathbb Z_2 \to  B_{r}(S^2)\to G(S^2_r)\to G(S^2)\to 1$$
\noindent where $S^2_r=S^2-r$  open disks.

It is well known that $G(S^2)=\mathbb Z_2$, where the nontrivial
 element is represented by
the isotopy class of an orientation reversing homeomorphism of the
 sphere. So the preimage of the element $id\in G(S^2)$ with respect to the
 projection  $G(S^2_r)\to G(S^2)$ gives
 the mappping class group  $Mod_{S^2_r}$(see section 4 for the definition
 of $Mod_S$).    The above  sequence implies   immediately the following
short exact  sequence:

 $$1\to \mathbb Z_2 \to  B_{r}(S^2)\to  Mod_{S^2_r} \to 1.$$

\noindent This sequence will  used in section 5.

\vskip 1.0cm

\section{ Automorphisms of the symplectic group }

Let $Sp(2n,\mathbb Z)$ denote the  symplectic group
  over the integers.  This group is sometimes   called   {\it proper
 symplectic group} as in  \cite{kms}. For some  choice of a basis,  the
 group
$Sp(2n,\mathbb Z)$  may  be defined as  the set of  matrices $2n\times
 2n$ over $\mathbb Z$ which satisfy  the following conditions: for all
 $1\leq l<j\leq 2n$,  
$\Sigma_{i=1}^{i=n}(a_{2i-1,l}a_{2i,j}- a_{2i-1,j}a_{2i,l})$ is equal  1 if
 $(l,j)$ is of the form $(2k-1,2k)$ and  is equal to 0 otherwise. In
 \cite{kms} and  \cite{ne}  the description of the symplectic group
  is somewhat  different to  the one given above. One definition is given in
 terms
 of a basis  and the other is given in terms of a basis  which corresponds
 to a permutation of the first base.
Recall that the elements of the the group $Sp(2n,\mathbb Z)$ are
 automorphisms which  are  obtained as the induced  homomorphisms
 in $H^1(S,\mathbb Z)$  by  an orientation preserving
 homeomorphisms of the
 orientable  closed  surface $S$ of genus $n$.

 We refer to \cite{ne}
and \cite{rein}
 for most of the properties of the group of symplectic  matrices. The
 study of  the Reidemeiser number of  automorphism of the  unimodular
 groups is   in progress, while  the corresponding  basic material may be 
 found  in
  \cite{hure} and \cite{ne}.

Following \cite{rein} let us call a homomorphism of  $Sp(2n,\mathbb Z)$
 into $\{\pm 1\}$ a {\it character}. It is known  that every
 automorphism $\phi$ of $Sp(2n,\mathbb Z)$ is given by
  $\phi(M)=\psi(M)UMU^{-1}$ for all
$M\in Sp(2n,\mathbb Z)$, where $\psi$ is a character and  $U$ is a
 matrix
 in  the full symplectic group.
We firstly  show that, if we change an automorphism by a character, then  we do
 not change the finiteness of the Reidemeister number  of the resulting
 automorphism. Thus   we reduce the problem to the same problem for the automorphisms given by conjugation
 by the elements of the full symplectic group. Since the proper
 sympletic group has index 2 in the full symplectic group,  then, by the  Lemma
 2.1,  we
 reduce our   problem to the  study of the finiteness of Reidemeister number of  ONE  automorphism
 given as  conjugation by
 one  element of the full symplectic group, which is not in the proper
 symplectic group.
We also need to know  the number of usual conjugacy classes in  the
group.

{\bf Lemma 3.1} Let $\phi:Sp(2n,\mathbb Z)\to Sp(2n,\mathbb Z)$ be an
 automorphism and  
$\psi:Sp(2n,\mathbb Z) \to \mathbb Z_2$ be a character.
$R(\phi)$ is finite if and only if $R(\psi\cdot\phi)$ is finite,
where homomorphism $\psi\cdot \phi$ is given by $(\psi\cdot \phi)(g)=\psi(g)\cdot\phi(g)$ for all  $g\in Sp(2n,\mathbb Z)$.
 \begin{proof} If $\psi$ is the trivial homomorphism then the result is
 clear since the two homomorphisms are the same. Otherwise, let $H$ be
 the kernel of $\psi$, which is a subgroup of index 2 in
 $Sp(2n,\mathbb Z)$. We split each Reidemeister class of $\phi$ into
 two disjoint
 subsets, namely two elements $y_1,y_2$ of a Reidemeister class are in
 the
 same subset if and only  if $y_2=ay_1\phi(a)^{-1}$ for some $a \in H$.
 It is
 straightforward to see that this relation gives  a partition of a
 Reidemeister class into two nonempty subsets. The set all these
 subsets as we
 run over the set of all Reidemeister classes of $\phi$ gives  a
 partition
 of the group. Now a Reidemeister class of the homomorphism
$\psi\cdot \phi$ which contains   an element $y$ is the union of the subset
 containing  $y$ with the
 subset containing  the element $ay\phi(a)^{-1}$ for some  element $a$
 which does not belong  to  $H$. This implies the result.
\end{proof}

Let $\phi$ be the automorphism which is the  conjugation by the
 diagonal matrix of order $2n$,  where the elements of the diagonal are
  $a_{i,i}=(-1)^{i+1}$. It is easy to verify that if  $M=(m_{i,j})$
 then
 $\phi(M)_{i,j}=(-1)^{i+j}m_{i,j}$.

Now we prove  the result for the case  $n=1$. The proof of the general case will 
follow the same strategy  as in  this case  and
will eventually    reduce   to the case $n=1$.

{\bf Proposition 3.2} If  $\phi:Sp(2,\mathbb Z)\to
 Sp(2,\mathbb Z)$ is an  automorphism, then
$R(\phi)=\infty$.
\begin{proof} Let $p$ be a prime. The set of invertible elements in the
 cyclic group
$\mathbb Z_p$ is a multiplicative group of  order $p-1$ denoted by
 $\mathbb Z_p^*$.
The number of elements of the set
$V_1=\{w\in \mathbb Z_p^*| w^2=- 1 \}$ is at most 2.
  If $p-1$ is divisible by 4 then   this set has 2 elements, otherwise it is empty .     Consider the groups $Sp(2,\mathbb Z_p)$,
 the
 Symplectic group module the prime $p$, when   $p$ runs over the set of
 all
 primes. The groups  $Sp(2,\mathbb Z_p)$ are   finite. We will use
 Lemma
 2.2 to obtain information about the automorphism $\phi$. The
 automorphism $\phi$   sends the
 matrix
$$
\begin{pmatrix}
a & b \\
c & d
\end{pmatrix}
$$
to
$$
\begin{pmatrix}
a & -b \\
-c & d 
\end{pmatrix}.
$$

The order of the  group $Sp(2,\mathbb Z_p)$ is  $p(p^2-1)$. This group  contains a
 subgroup isomorphic to $Z_{p-1}$. Consider  the subgroup generated by the
 diagonal matrices of the form

$$
\begin{pmatrix}
w & 0 \\
0 & w^{-1} 
\end{pmatrix}
$$

\noindent where $0\ne w\in \mathbb Z_p$. Denote by $\bar w$ the above  matrix
  determined by the element $w$. We  now compute  the    elements
 of this subgroup
 which are  in the same  Reidemeister class of $\bar w$ for $w$ not in
 $V_1$. We  claim that  $\bar w$ and $-\bar w^{-1}$ are the only elements
 in the Reidemeister class of 
$\bar w$. This will
 imply  that the number of Reidemeister classes of $\phi$ is at least
 $(p-3)/2$. For the calculation of the elements of the subgroup which are
 in
 the Reidemeister class of $\bar w$ let us take an arbitrary matrix $M$
 in
 $Sp(2,\mathbb Z_p)$ and let us require  that
$ M\bar w\phi(M^{-1})$ is diagonal. If the Symplectic matrix $M$ is
 given by 

$$
\begin{pmatrix}
a & b \\
c & d 
\end{pmatrix}
$$

\noindent then $\phi(M)$ is given by 

$$
\begin{pmatrix}
a & -b \\
-c & d 
\end{pmatrix}
$$

\noindent and $\phi(M^{-1})$ is equal to 

$$
\begin{pmatrix}
d & b \\
c & a 
\end{pmatrix}
$$

\noindent and the product $ M\bar w\phi(M^{-1})$ is
equal to

$$
\begin{pmatrix}
wad+w^{-1}bc & wab+w^{-1}ab \\
wcd+w^{-1}cd & wbc+w^{-1}ad 
\end{pmatrix}.
$$

 The above matrix is diagonal if and only if   $(w+w^{-1})ab=0$
 and 
$(w+w^{-1})cd=0$. Since $w+w^{-1}\ne 0$  because $w$
 does not belong to  $V_1$, we must have $ab=0$ and $cd=0$. Since  $M$ has
 determinant 1,
the possible solutions are: 1) $a=0$ and $d=0$  or 2) $b=0$ and $c=0$.
 In the first case
 we  have $bc=-1$ and the diagonal  elements   are
 $-w^{-1}$ and $-w$. 
 In the second case  we have $ab=1$ and the diagonal  elements  are $w$ and $w^{-1}$. 
So each Reidemeister class
 contains  2 elements of this form. So the number of Reidemeister
 classes
 is at most $(p-3)/2$ and the result follows.

\end{proof}

Now we prove  the main result of this section.

{\bf Theorem 3.3}  If  $\phi:Sp(2n,\mathbb Z)\to
 Sp(2n,\mathbb Z)$ is an  automorphism, then  
$R(\phi)=\infty$.

\begin{proof} We will reduce the problem to the case of $Sp(2,\mathbb
 Z)$ at the very  end. 
Let $M \in Sp(2n,\mathbb Z)$ be given by

$$
\begin{pmatrix}
a_{1,1} & a_{1,2} &.....& a_{1,2n-1} & a_{1,2n} \\
b_{1,1} & b_{1,2} &.....& b_{1,2n-1} & b_{1,2n} \\
. & . & .   &.  &  \\
. & . & .   &.  &  \\
a_{n,1} & a_{n,2} &.....& a_{n,2n-1} & a_{n,2n} \\
b_{n,1} & b_{n,2} &.....& b_{n,2n-1} & b_{n,2n} 
\end{pmatrix}.
$$

Then $\phi(M)^{-1}$ is given by 

$$
\begin{pmatrix}
b_{1,2} & a_{1,2} &.....& b_{n,2} & a_{n,2} \\
b_{1,1} & a_{1,1} &.....& b_{n,1} & a_{n,1} \\
. & . & .   &.  &  \\
. & . & .   &.  &  \\
b_{1,2n} & a_{1,2n} &.....& b_{n,2n} & a_{n,2n} \\
b_{1,2n-1} & a_{1,2n-1} &.....& b_{n,2n-1} & a_{n,2n-1}
\end{pmatrix}.
$$

Let $\bar w$ be the matrix

$$
\begin{pmatrix}
w& 0 &.....& 0 & 0 \\
0 & w^{-1}& ...& 0&0\\
. & . & .   &.  &  \\
 . & .   &.  &  \\
0 & 0 &.....& 1 & 0 \\
0 & 0 &.....& 0 & 1
\end{pmatrix}.
$$
 In order to describe $ M\bar w\phi(M^{-1})$,  let us introduce the
 following notation.
 Let
$$A_{2i-1,2j-1}=w^{\epsilon}a_{i,1}b_{j,2}+
 w^{-\epsilon}b_{j,1}a_{i,2}..+..a_{i,2n-1}b_{j,2n}+a_{i,2n}b_{j,2n-1},$$
$$A_{2i-1,2j}=w^{\epsilon}a_{i,2}a_{i,1}+w^{-\epsilon}a_{i,2}a_{i,1}..+..a_{i,2n-1}a_{i,2n}+a_{i,2n}a_{i,2n-1},$$
 $$ A_{2i,2j-1}=w^{\epsilon}b_{j,1}b_{j,2}+
 w^{-\epsilon}b_{j,1}b_{j,2}.+..b_{j,2n-1}b_{j,2n}+b_{j,2n}b_{j,2n-1},$$
$$ A_{2i,2j}=
 wb_{j,1}a_{i,2}+w^{-1}b_{j,2}a_{i,1}..+..b_{j,2n-1}a_{i,2n}+b_{j,2n}a_{i,2n-1}.$$
\noindent where $\epsilon=1$  for $j=1,2$ and zero otherwise.\\

The product $ M\bar w\phi(M^{-1})$ is given by
\smallskip
$$
\begin{pmatrix}

A_{1,1}
 &
A_{1,2}
& .....& A_{1,2n-1}
&
A_{1,2n}
\\
A_{2,1}

&
A_{2,2}

&.....&
A_{2,2n-1}
&
A_{2,2n}
 \\
...&.. . &.. .   &. .. &  \\
... & ... & . ..  &...&  \\
A_{2n-1,1}
&
A_{2n-1,2}
&.....& A_{2n-1,2n-1}
&
A_{2n-1,2n}
\\
A_{2n,1}
&
A_{2n,2}
&.....& A_{2n,2n-1}
&
A_{2n,2n}
\end{pmatrix}.
$$
\noindent

\smallskip

It  follows that $wa_{i,1}b_{j,2}+
w^{-1}b_{j,1}a_{i,2}..+..a_{i,2l-1}b_{j,2l}+a_{i,2l}b_{j,2l-1}..+..a_{i,2n-1}b_{j,2n}+a_{i,2n}b_{j,2n-1} = wa_{i,2}b_{j,1}+
 w^{-1}b_{j,2}a_{i,1}..+..a_{i,2l}b_{j,2l-1}+a_{i,2l-1}b_{1,2l}..+.. 
a_{i,2n}b_{j,2n-1}+a_{i,2n-1}b_{j,2n}=0$  for all $(i,j)\ne (1.1)$,
 $1\leq i,j \leq n$.
This implies the system of equations:
$(w-w^{-1})b_{j,2}a_{i,1}=(w-w^{-1})b_{j,1}a_{i,2}$ for  all $(i,j)\ne
 (1,1)$, $1\leq i,j \leq n$. So all $2\times 2$ submatrices of the
 matrix whose columns
are $a_{1,1}, b_{1,1},....a_{n,1}, b_{n,1}$ for $1\leq i\leq n$
and  $a_{1,2}, b_{1,2},....a_{n,2}, b_{n,2}$ for $1\leq i\leq n$
 $b_{i,1}$ for $1\leq i\leq n$ different from

$$
\begin{pmatrix}
a_{1,1} & b_{1,1}  \\
a_{2,1} & b_{2,1}
\end{pmatrix}
$$
have determinant zero.
This implies that   the above  $2\times 2$  matrix
   has determinant 1( so is  not zero) and
  $a_{i,1}=a_{i,2}=b_{i,1}=b_{i,2}=0$ for $i>1$.

So the  matrices $M$ and $\bar w\phi(M)^{-1}$ are respectively  of the form

$$
\begin{pmatrix}
a_{1,1} & a_{1,2} & a_{1,3} & a_{1,4} &.....& a_{1,2n-1} & a_{1,2n} \\
b_{1,1} & b_{1,2} & b_{1,3} & b_{1,4} &.....& b_{1,2n-1} & b_{1,2n} \\
0       & 0       & a_{2,3} & a_{2,4} &.....& a_{2,2n-1} & a_{2,2n} \\
 0      & 0       & b_{2,3} & b_{2,4} &.....& b_{2,2n-1} & b_{2,2n} \\

.       & .             & .   &.  &.... .& . &. \\

 0 & 0 & a_{n,3} & a_{n,4} & .....& a_{n,2n-1} & a_{n,2n} \\
0 & 0 & b_{n,3} & b_{n,4} & .....& b_{n,2n-1} & b_{n,2n}
\end{pmatrix}.
$$

and 

$$
\begin{pmatrix}
wb_{1,2} & w^{-1}a_{1,2} & 0 &0 &.....& 0 & 0 \\
wb_{1,1} & w^{-1}a_{1,1} & 0& 0& .....& 0 & 0 \\
b_{1,4} & a_{1,4} & b_{2,4} & a_{2,4} &.....& b_{n,4} & a_{n,4} \\
b_{1,3} & a_{1,3} & b_{2,3} & a_{23} &.....& b_{n,3} & a_{n,3}  \\
. & . & .   & .  & ..... & . & .  \\
. & . & .   &.  &..... & . & .  \\
b_{1,2n} & a_{1,2n} &b_{2,2n} & a_{2,2n} &.....& b_{n,2n} & a_{n,2n} \\
b_{1,2n-1} & a_{1,2n-1} &b_{2,2n-1} & a_{2,2n-1} &.....& b_{n,2n-1} &
 a_{n,2n-1} 
\end{pmatrix}.
$$

 The product $M\bar w(\phi(M)^ {-1}$ is  of the form

$$
\begin{pmatrix}
A & 0  \\
0 & I_{2n-2}  
\end{pmatrix}
$$

\noindent where the A is of order $2\times 2$, $I_{2n-2}$ is the
 identity matrix of order $2n-2$, and $0's$ are  the  trivial matrices
 of
 orders $2\times 2n-2$, $2n-2\times 2$, respectively. Hence  the submatrix
 of
 order $2n-2\times 2n-2$

$$
\begin{pmatrix}
 b_{2,4} & a_{2,4} &.....& b_{n,4} & a_{n,4} \\
 b_{2,3} & a_{4,3} &.....& b_{n,3} & a_{n,3}  \\
 .   & .  & ..... & . & .  \\
.   &.  &..... & . & .  \\
b_{2,2n} & a_{2,2n} &.....& b_{n,2n} & a_{n,2n} \\
b_{2,2n-1} & a_{2,2n-1} &.....& b_{n,2n-1} & a_{n,2n-1} 
\end{pmatrix}
$$

\noindent is invertible. Regarding  the columns as a vector in the 2n-2
  dimensional vector space over the rationals,   these 2n-2
  vectors form a basis. If $v$ is any of the above column, then the
 inner product
of $(a_{1,1} , a_{1,2} , a_{1,3} , a_{1,4} ,....., a_{1,2n-1} ,
 a_{1,2n})$ with the column
 vector $(0,0, v)$
 is zero. Therefore the inner product of $( a_{1,3} , a_{1,4} ,.....,
 a_{1,2n-1} , a_{1,2n})$ and $v$  is also zero for all $v$. Since  the
 set
 of all $v's$ forms a basis, this implies that $( a_{1,3} , a_{1,4}
 ,.....,
 a_{1,2n-1} , a_{1,2n})$ is trivial and the matrices $M$ and $\bar
 \phi(M)^{-1}$ are  of the form 
$$
\begin{pmatrix}
a_{1,1} & a_{1,2} & 0& 0&.....& 0& 0\\
b_{1,1} & b_{1,2} & 0       & 0       &.....& 0          & 0        \\
0       & 0       & a_{2,3} & a_{2,4} &.....& a_{2,2n-1} & a_{2,2n} \\
 0      & 0       & b_{2,3} & b_{2,4} &.....& b_{2,2n-1} & b_{2,2n} \\

.       & .             & .   &.  &.... .& . &. \\

 0 & 0 & a_{n,3} & a_{n,4} & .....& a_{n,2n-1} & a_{n,2n} \\
0 & 0 & b_{n,3} & b_{n,4} & .....& b_{n,2n-1} & b_{n,2n}
\end{pmatrix}.
$$

and 

$$
\begin{pmatrix}
wb_{1,2} & w^{-1}a_{1,2} & 0 &0 &.....& 0 & 0 \\
wb_{1,1} & w^{-1}a_{1,1} & 0& 0& .....& 0 & 0 \\
0       & 0       & b_{2,4} & a_{2,4} &.....& b_{n,4} & a_{n,4} \\
0       & 0       & b_{2,3} & a_{23} &.....& b_{n,3} & a_{n,3}  \\
. & . & .   & .  & ..... & . & .  \\
. & . & .   &.  &..... & . & .  \\
0        & 0        &b_{2,2n} & a_{2,2n} &.....& b_{n,2n} & a_{n,2n} \\
0          & 0          &b_{2,2n-1} & a_{2,2n-1} &.....& b_{n,2n-1} &
 a_{n,2n-1} 
\end{pmatrix}.
$$

 The $2\times 2$  submatrix $A$  of
the product $M\bar w(\phi(M)^ {-1}$   is given  as in the
 case when  $n=1$.
Thus   we will have as solution   for $a_{1,1},a_{1,2}, b_{1,1}, b_{1,2}$
 the same solutions
as in the case $n=1$ and  the result follows. 
  
\end{proof}

\vskip 1.0cm

\section {Automorphisms of the mapping class group}

\vskip 1.0cm
 In this section we study the $R_{\infty}$ property for  the  Mapping class group of a closed surface. 
As an application of the
results of the  previous section, we show that  mapping class group of a closed orientable surface of
 genus $g>0$ has the $R_{\infty}$ property.
 We start by quoting some   results about Outer automorphism group  of the mapping class group.
 Let  $S$ be a compact orientable surface. Recall that
 $Mod_S$  denotes  the group of orientation-preserving
  homeomorphisms of $S$ modulo isotopy,
and  $Mod_S^*$  denotes the group of the homeomorphisms of $S$ modulo
 isotopies. From \cite{iva} we have

{\bf Theorem 4.1} If $S$ is neither  a sphere with $\leq 4$ holes, nor
 a
 torus with $\leq 2$ holes, nor a closed surface of genus 2, then every
 automorphism of $Mod_S$ is given by the restriction of an inner
 automorphism of $Mod_S^*$. In particular $Out(M_S)$ is a finite group
 and
 moreover,
$$Out(Mod_S)=Z_2, Out(Mod_S^*)=1.$$
 
{\bf Theorem 4.2} If $S$ is   a closed surface of genus 2  then
 $Out(M_S)$ is canonically isomorphic to $H^1(Mod_S;Z_2)\oplus
 Mod_S^*/Mod_S=Z_2\oplus Z_2$ and 
$Out(M_S^*)$ is canonically isomorphic to $H^1(Mod_S^*;Z_2)=Z_2\oplus
 Z_2$.

   In  \cite{iva} the remaining  cases  also are given. In the cases where
 $S$ is either the torus $T$ or the torus with one hole  we have
 $Mod_S=Sl_2(\mathbb Z)$ and  $Mod_S^*=Gl_2(\mathbb Z)$,
 so $Out(Mod_S^*)=1$ and $ Out(Mod_S)=Z_2.$ 
For $S=S^2$ we have  $Mod_S=\{1\}$, the trivial group,  and 
$Mod_S^*=\mathbb Z_2$, therefore
$Out(Mod_S^*)=Out(Mod_S)=1.$ 
 
We will show that the mapping class group of closed surface has the  $R_{\infty}$ property.

{\bf Theorem 4.3} Let $S$ be an orientable  closed surface. The mapping
 class group $Mod_S$ has the $R_{\infty}$ property if and only if $S$
 is not $S^2$.

\begin{proof}  If $S$ is $S^2$  we know that the mapping class group is
  finite. If $S=T$ then
 $Mod_S=Sl_2(\mathbb Z)$ and this group   has the
 $R_{\infty}$ property. In the remaining  cases,  we consider the short exact
 sequence
$$1\to I_n \to Mod_S \to Sp(2n,\mathbb Z)\to 1,$$

\noindent where $I_n$ is the kernel of the  well- known   homomorphism
 $p: Mod_S \to Sp(2n,\mathbb Z)$, also known  as the Torelli group.

If the genus  $g>2$,  then we consider an automorphism  $\phi$ which represents
 the
 nontrivial element of $Out(Mod_S)=Z_2$. Pick  an  automorphism $\phi$ which is  
  conjugation by an  element of the group  $Mod_S^*$  represented by a  orientation reversing homeomorphism of surface $S$. It preserves the Torelli group, and so we obtain   a homomorphism   of the   short exact sequence. 
Thus we obtain an  induced automorphism $ \ovl{\phi}$ in the quotient group $Sp(2n,\mathbb Z)$.

 The Reidemeister number  $R(\ovl\phi)$ is infinite by Theorem 3.3 . Then remark about extensions (\ref{extens})  implies that the Reidemeister number $R(\phi)$ is also infinite. Then the $R_{\infty}$ property   follows from Lemma 2.1.
For genus $g= 2$, we have to show that 3  automorphisms of $Mod_S$ have infinite numbers of twisted
conjugacy classes,
 since $Out(ModS)=\mathbb Z_2+\mathbb Z_2$ in this case. These three automorphisms are all
 conjugation, so  they   preserve the Torelli group and  the result
 follows as above.

\end{proof}

\section{Automorphisms of the  braid groups of  $S^2$ and disk  $D$}
 
 We  prove in  this section the   $R_{\infty}$ property for $B_n(S^2)$.
 This case will be
 investigated  using the results of the previous section  and
 the relation between
 mapping class groups and braids groups  given in  section 2. 
 Then we will consider the
 case when  $S$ is a disk $D$, namely  the group  $B_n(D)$,
 also called   the
 Artin braid group on  $n-$strings.

{\bf Theorem 5.1} The group  $B_n(S^2)$ has  the $R_{\infty}$ property
 if and only if $n>3$.

\begin{proof} Let $\phi: B_n(S^2) \to B_n(S^2)$ be an automorphism.
 Since the center of   $B_n(S^2)$  is a characteristic subgroup, $\phi$ induces a homomorphism of the short  exact sequence 

$$1\to \mathbb Z_2 \to  B_{n}(S^2)\to Mod_{S^2_n} \to 1 $$
\noindent
 where the  short exact sequence was   obtained from the sequence
 in Theorem 2.3. The results of  Appendix  imply
 that  the group $Mod_{S^2_n}$,  for  $n>3$,  has 
the $R_{\infty}$ property. Then the remark  about extensions (\ref{extens}) implies that  the group  $B_{n}(S^2)$ also has this property.
 For $n\leq 3$ the groups $B_n(S^2)$ are finite so they do not have the
 $R_{\infty}$ property.
\end{proof}

{\bf Remark} For $n=4$ the lower central series of $B_4(S^2)$ has been
 studied. The commutator subgroup $[B_4(S^2), B_4(S^2)]$ is isomorphic
 to 
$Q_8\rtimes F_2(x,y)$, the semi-direct product of the quatenionic group
 with the free group on two generators \cite{gg}, whose Abelinization  is $\mathbb
 Z_6$. Since the  $[B_4(S^2), B_4(S^2)]$ 
is characteristic, to show that $B_4(S^2)$ has the $R_{\infty}$
 property, it is sufficient  to show that  $[B_4(S^2), B_4(S^2)]\simeq Q_8\rtimes F_2(x,y)$ has $R_{\infty}$ property. But, using the remark  about
extensions (\ref{extens}),  this follows  from the fact that
$Q_8 \subset  Q_8\rtimes F_2(x,y)$ is characteristic and  that group
 $F_2(x,y)$ has the $R_{\infty}$ property.

Now let us consider the Artin braid group $B_n(D)$. We will show that 
$B_n(D)$ has the $R_{\infty}$ property, except a few cases.  It is known
  from \cite{dy-gr}
  that $Out(B_n(D))=\mathbb Z_2$, but  we  do not   use this fact.

 Denote by $\mathbb B_n(D)$ the quotient of $B_n(D)$ by the center.

{\bf Proposition 5.2} The group  $\mathbb B_n(D)$   has the
 $R_{\infty}$
 property if and only if  $n>2$.
\begin{proof} The ``if part" follows from the fact that $B_1(D)$ is trivial and $B_2(D)$ is   isomorphic to
 $\mathbb Z$. So let $n>2$. From Theorem 15 in
 \cite{dy-gr} we know that the free group $F_{n-1}$ of rank 
$n-1$ is characteristic in $\mathbb B_n(D)$. From the formulas which
 relates the Reidemeister classes  of the terms of a short exact
 sequence(see \cite{go:nil1})
 we conclude that    $\mathbb B_n(D)$ has the $R_{\infty}$ property.
\end{proof}

{\bf Theorem 5.3} The group  $B_n(D)$ has the $R_{\infty}$ property if
 and
 only if  $n>2$.
\begin{proof} The ``if part" follows from the fact that $B_1(D)$ is trivial and $B_2(D)$ is   isomorphic to
 $\mathbb Z$. So let $n>2$.  Consider the short
 exact sequence

$$1\to \mathbb Z \to  B_n(D) \to \mathbb B_n\to 1 .$$

\noindent
 Because $\mathbb Z$ is the center of  $B_n(D)$ and the
 center is characteristic, any automorphism of $B_n(D)$ is an
 automorphism
 of short exact sequence. Since $\mathbb B_n$ has the $R_{\infty}$
 property the result follows from  remark  about
extensions (\ref{extens}).
\end{proof}

\section{ Appendix: Geometric group theory  and $R_{\infty}$ property for mapping class group}

The object of this Appendix  is to use geometric  methods to  prove more stronger result then in section 4 .
 Namely, we show  that  the  known geometric
methods apply to prove that  for all compact orientable  surfaces 
with genus $g$ and $p$ boundary components, and $3g+p-4>0$  
the mapping class group  has the $R_{\infty}$ property with a  few exceptions.

 We need to use the non-elementary 
result of Masur and Minsky
\cite{MM} (see also Bowditch \cite{B}) that the complex of curves of an
oriented surface, with genus $g$ and $p$ boundary components, and
$3g+p-4
>0$,  is Gromov-hyperbolic space.

{\bf Lemma 6.1 }
If $G$ is a group and $\varphi$
 is an endomorphism of $G$,  then
 an element $x\in G$ is always
 $\varphi$-conjugate to its image $\varphi(x)$.
\begin{proof} 
Put $\gamma=x^{-1}$. Now $x$ is $\varphi$-conjugate to
 $x^{-1} x \varphi(x) =\varphi(x)$.
\end{proof}

Let $G$ be a group, and $\varphi$ an automorphism of $G$
of order $m$. Construct  the group
$G_\varphi= G\rtimes_\varphi \Z_m = \langle G,t | \, \forall g\in G,
tgt^{-1} = \varphi(g),\, t^m=1 \rangle $.

{\bf Lemma 6.2 }
Two elements $x,y$ of $G$ are $\varphi$-conjugate if and only if  $xt$
 and $yt$ are conjugate
in the usual sense in $G_\varphi$. Therefore $R(\varphi)$ is the number
 of usual conjugacy
classes in the coset $G\cdot t$ of $G$ in $G_\varphi$.
\begin{proof} 
If $x$ and $y$ are $\phi$-conjugate,  then there is a $\gamma \in G$
 such
that $\gamma x=y\varphi(\gamma)$. This implies $\gamma x=yt\gamma
 t^{-1}$
and therefore $\gamma(xt)=(yt)\gamma$.   So $xt$ and $yt$ are conjugate
 in the usual sense in $G_\varphi$.
Conversely suppose $xt$ and $yt$ are conjugate in $G_\varphi$.
Then there is a $\gamma t^n \in G_\varphi$ with $\gamma t^n xt=yt\gamma
 t^n$.
From the relation $txt^{-1}=\varphi(x)$,  we
obtain $\gamma \varphi^n(x)t^{n+1}=y\varphi(\gamma) t^{n+1} $
 and therefore $\gamma \varphi^n(x)=y\varphi(\gamma)$.
This shows that $\varphi^n(x)$ and $y$ are $\varphi$-conjugate.
However,  by lemma 6.1 ,   $x$ and $\varphi^n(x)$ are
 $\varphi$-conjugate, so $x$ and $y$
must be $\varphi$-conjugate.

\end{proof}

The following lemma was proven by Delzant

{\bf Lemma 6.3 } \cite[Lemma 3.4]{ll}
If $K$  is a normal subgroup of a group $\Gamma$
acting
non-elementary on a hyperbolic space, and if $\Gamma/K$ is abelian,
then
any coset of $K$ contains infinitely many conjugacy classes.

\begin{proof} 
Fix $u$ in the coset $C$ under consideration. Suppose for a moment that
we can find $c,d \in K$, generating a free group of rank 2, such that
$uc^{\infty}\not=c^{-\infty}$ and $ud^{\infty}\not=d^{-\infty}$
(recall that we denote $g^{-\infty}=\lim_{n\rightarrow +\infty}g^{-n}$
for $g$ of infinite order). Consider $x_k=c^kuc^k$ and $y_k=d^kud^k$.
For $k$ large, the above inequalities imply that these two elements
 have infinite order, and do not generate a virtually cyclic group because
 ${x_k}^{+\infty}$ and ${x_k}^{-\infty}$(respectively ${y_k}^{+\infty}$
 and ${y_k}^{-\infty}$ ) is close  to $c^{+\infty}$ and $c^{-\infty}$ (
respectively $d^{+\infty}$ and $d^{-\infty}$). Fix $k$, and consider
the elements  $z_n={x_k}^{n+1}{y_k}^{-n}$. They belong to the coset
 $C$,
because $\Gamma/K$ is abelian, and their stable norm goes to infinity
 with $n$.
Therefore $C$ contains infinitely many
conjugacy classes.
 Let us now construct $c,d$ as above. Choose $a,b \in K$ generating a
 free group
of rank 2. We first explain how to get $c$. There is a problem only
if $ua^{\infty}=a^{-\infty}$ and $ub^{\infty}=b^{-\infty}$ . In that
 case
there exists integers $p,q$ with $ua^pu^{-1}=a^{-p}$ and
  $ub^qu^{-1}=b^{-q}$.
 We take $c=a^pb^q$, noting that $ucu^{-1}=a^{-p}b^{-q}$ is different
 from
$c^{-1}=b^{-q}a^{-p}$.

Once we have $c$, we choose $c^*\in K$ with $<c,c^*>$ free of rank 2,
 and we
obtain $d$ by applying the preceding argument using $c^*$ and $cc^*$
 instead
of $a$ and $b$. The group $<c,d>$ is free of rank 2 because $d$ is a
 positive
word in $c^*$ and $cc^*$.

\end{proof}

{\bf Theorem 6.4}
 If $G_\varphi$ has a non-elementary action by isometries on a
Gromov-hyperbolic length space, then $G$ has infinitely many
$\varphi$-twisted conjugacy classes.

\begin{proof} 
By elementary action we mean an action consisting of elliptic elements,
or
with a global fixed point, or with a global fixed pair  in the boundary of
the
hyperbolic space.
The  statement of the theorem immediately follows from Lemma 6.2
and Delzant Lemma 6.3.
\end{proof}

The Theorem 6.4
applies if $G$ is a Gromov-hyperbolic group and
 $\varphi$ has finite order in $Out(G)$. In fact in this case,
$G_\varphi$ contains $G$ as a subgroup of finite index, thus is
quasi-isometric to $G$, and by quasi-isometry invariance, it is itself
Gromov-hyperbolic group. 

The Theorem 6.4
also applies when $G$ is so-called relatively
hyperbolic and
$\varphi$ is of finite order in $Out(G)$. In this case, the
quasi-isometry
invariance   of $G_\varphi$ and $G$  is harder to establish, but is has been
 proven by Dru\c{t}u in
\cite{Dr}.
 
{\it Remark. } For automorphism of a hyperbolic group, of infinite
order in
$Out(G)$, the fact that  $G$ has infinitely many $\varphi$-twisted
conjugacy
classes was implicitly  proven by Levitt and Lustig in \cite{ll},
see also \cite{FelPOMI}. Similarly we can prove that $R(\varphi)=\infty$
for automorphism $\varphi$ of infinite order in $Out(G)$, where $G$ is 
relatively hyperbolic group.

Let now $S$ be an oriented compact surface with genus $g$ and with  $p$
boundary components, where  $3g+p-4 >0$. It is easy to see that the
 mapping class group
$Mod_S$ is a normal subgroup of  the full mapping class group $
 Mod_S^*$, of  index  $2$.
The graph of curves of $S$,  denoted $\mathcal{G}(S)$,  is the
graph
whose vertices are the simple curves of $S$ modulo isotopy. Two
vertices (that is two isotopy classes of simple curves) are linked by
an
edge in this graph if they can be realized by disjoint curves. Both
$Mod_S$ and  $ Mod_S^*$ act on $\mathcal{G}(S)$ in a
non-elementary way.

Thus, Theorem 6.4 
is applicable  for $Mod_S$ and
 for $\varphi_1$  the automorphism induced by reversing the
orientation of $S$, since in this  case,
$({Mod_S})_{\varphi_1}={{Mod_S}\rtimes}_{\varphi_1} \Z_2 \simeq
 Mod_S^*$.
For $Mod_S$ and
$\varphi_0= Id$ we have $R(\varphi_0= Id)=\infty$ because the group
  $Mod_S$
has infinite number of usual conjugacy classes.
Finally,
 $Out(Mod_S) \simeq \{\overline{\varphi_0}, \overline{\varphi_1}
\}$ see  \cite{iva},
which  by Lemma 2.1 ensures   that $Mod_S$ has property $R_\infty$
if  $S$ is  an orientable  compact surface of  genus $g$ with  $p$
boundary components, where  $3g+p-4 >0$. The only cases not covered by
 this inequality are where  S is the torus with at most one hole  or where  S is  the sphere with at most 4
 holes. The case of a torus with at most one hole  follows from
the  section 3. The case of the  sphere with at most 4 holes   follows directly
 from the knowledge of Out(ModS)  and the cardinality of the mapping class group.

\end{document}